# A QUARTIC SURFACE OF INTEGER HEXAHEDRA

ROGER C. ALPERIN

ABSTRACT. We prove that there are infinitely many six sided polyhedra in $\mathbf{R}^3$, each with four congruent trapezoidal faces and two congruent rectangular faces, so that the faces have integer sides and diagonals, and also the solid has integer length diagonals. The solutions are obtained from the integer points on a particular quartic surface.

A long standing unsolved problem asks whether or not there can be a parallelipiped in $\mathbf{R}^3$ whose sides and diagonals have integer length. If one weakens the requirement and just asks for a six-sided polyhedron with quadrilateral faces, then one can find examples with integer length sides and diagonals. Peterson and Jordan [1], described a method for making these 'perfect' hexahedra. We review their method.

One takes two congruent rectangles and places them in space so that they are parallel, with the bottom rectangle rotated ninety degrees from the position of the top rectangle. Now connect the sides of the two rectangles with four congruent trapezoids. The shape is a piecewise linear version of the placement of two cupped hands together, at ninety degrees in clapping position. The centers of the rectangles lie on a line perpendicular to the top and bottom faces. If the sides of the rectangle have lengths $a, b$ then the diagonal has length $c$, where $a^2 + b^2 = c^2$. The parallel sides of the trapezoids are then also $a, b$. The slant side of the trapezoid is say $e$ and its diagonal is $d$. It follows from Ptolemy's Theorem that $d^2 = e^2 + ab$. Consider the trapezoid with base on the top rectangle of side $a$ and other base on the bottom rectangle opposite that edge of side $b$, having the slant sides of length $d$. Its diagonal is of length $f$, and also it is the interior diagonal of the hexahedron; thus $f^2 = d^2 + ab$. We shall refer to these as perfect hexahedra.

**Proposition 1.** *The simultaneous positive integer solutions to $a^2 + b^2 = c^2, d^2 = e^2 + ab, f^2 = d^2 + ab$ give the edge and diagonal lengths of a perfect hexahedron with two rectangular congruent opposite parallel faces, and four congruent trapezoidal faces.*

Peterson and Jordan asked if there are infinitely many such perfect hexahedra. They also gave several examples, including the small example $a = 8, b = 15, c = 17, e = 7, d = 13, f = 17$ and asked if it is the smallest. We provide a measure of size of solutions, and answer both of these questions affirmatively.

We can rewrite this set of equations. Basically, $d^2 - e^2 = f^2 - d^2 = ab$, so that we have three integer squares in arithmetic progression. This is the same as $e^2 + f^2 = 2d^2$ and $ab = f^2 - d^2$. Now then we have two norm equations $a^2 + b^2 = c^2, e^2 + f^2 = 2d^2$ over $Z[i]$.

It is well known that the relatively prime or primitive solutions to the Pythagorean equation, $x^2 + y^2 = z^2$ are given by $x = m^2 - n^2, y = 2mn, z = m^2 + n^2$ for (relatively prime) integers $m, n$. All integer solutions are scalar (integer) multiples of the primitve solutions. To obtain solutions to $x^2 + y^2 = 2z^2$, we can use solutions obtained from the solutions to the Pythagorean equation. Namely, corresponding to a solution form the complex number, $x + iy$, $|x + iy|^2 = x^2 + y^2$, and $|(x + iy)(1 + i)|^2 = 2(x^2 + y^2)$; thus a solution to the Pythagorean equation, yields after multiplication by $1 + i$ a solution to the second equation. Also a solution to $x^2 + y^2 = 2z^2$, gives the complex number, $x + iy$, so that $|(x + iy)\frac{(1-i)}{2}|^2 = \frac{x^2 + y^2}{2}$; thus the complex number $(x + iy)\frac{(1-i)}{2}$, provides the solution to the Pythagorean equation. This gives a bijection of the sets of solutions of these two equations. (Using this together with sign changes or interchanging of variables, accounts for all solutions.)

For the perfect hexahedron then, we parameterize the integer solutions to the equation $a^2 + b^2 = c^2$ as $a = \lambda(r^2 - s^2), b = \lambda(2rs)$. For the solutions to $e^2 + f^2 = 2d^2$ we parameterize with $p, q$, $e + if = \mu((p^2 - q^2) + i2pq)(1 + i)$, so that $e = \mu(p^2 - q^2 - 2pq), f = \mu(p^2 - q^2 + 2pq)$. The condition for a perfect hexahedron is that $ab = \frac{f^2 - e^2}{2}$.

**Proposition 2.** *Using these parameterizations, a perfect hexahedron is obtained from any integer solution to $2(p^2 - q^2)pq\mu^2 = (r^2 - s^2)rs\lambda^2$ with $r \neq \pm s, p \neq \pm q$, $\mu, \lambda, r, s, p, q \neq 0$ and conversely.*

**Proof:** As we have seen the perfect hexahedron gives rise to the equation $ab = \frac{f^2 - e^2}{2}$, which by our parameterization is a multiple of the equation $2(p^2 - q^2)pq\mu^2 = (r^2 - s^2)rs\lambda^2$. Conversely, given any non-trivial integer solution to this equation, we can form $a = \lambda(r^2 - s^2), b = 2rs\lambda, c = \lambda(r^2 + s^2), e = \mu(p^2 - q^2 - 2pq), f = \mu(p^2 - q^2 + 2pq), d = \mu(p^2 + q^2)$, which give non-zero integer solutions to the perfect hexahedron. ∎

First we shall determine the smallest solution. We measure the *size* of a solution by the number $\frac{ab}{2}$, which is the same as $|(r^2 - s^2)rs\lambda^2|$.

**Lemma 1.** *The size of any solution of the hexahedra equations is divisible by 60.*

**Proof:** Consider the equation $a^2 + b^2 = c^2$. Modulo 3 the squares are 0 or 1 so it is impossible that both $a^2$ and $b^2$ are both 1 $mod$ 3. Thus $ab$ is divisible by 3.

Consider this same equation $mod$ 5. Modulo 5 the squares are 0, 1 or 4. Therefore, the only solutions $mod$ 5 are $0 + 0 = 0$ or $1 + 4 = 0$; in either case, $abc$ is divisible by 5. Suppose if posssible that $ab$ is not divisible by 5. Using the parameterization described above, then 5 divides $c = \lambda(r^2 + s^2)$, but not any of $\lambda, r, s, r-s, r+s$. Also, $\frac{ab}{2} = 2pq(p^2 - q^2)\mu^2$, for a non-trivial hexahedron. So therefore 5 does not divide any of $p, q, p+q, p-q, \mu$. However, then the numbers $r, s, r+s, r-s$ are all different $mod$ 5 and nonzero; similarly, for $p, q, p+q, p-q$. Thus, for a solution to the perfect hexahedron, $48\mu^2 = 2(p^2-q^2)pq\mu^2 = (r^2-s^2)rs\lambda^2 = 24\lambda^2 \ mod$ 5; this is impossible. Therefore $ab$ is divisible by 5.

Finally, the size is $|2pq(p^2 - q^2)\mu^2|$ which is certainly even, but also if $p, q$ are both odd, $p^2 - q^2$ is even; so the size is always divisible by 4. ∎

**Theorem 1.** *There is a unique hexahedron with solution of size 60; it has $a = 8, b = 15, c = 17, d = 13, e = 7, f = 17$. The second smallest is a hexahedron with solution of size 120; it has $a = 24, b = 10, c = 26, d = 16, e = 7, f = 23$.*

**Proof:** Solving the equation $pq(p^2 - q^2)\lambda^2 = 30$ in positive integers, we see immediately that $\lambda = \pm 1$, and $p, q$ are divisors of 30; say $p > q$ and $p + q > p > q > p - q > 0$ (or possibly $p + q > p > p - q > q > 0$). Hence $p - q = 1$. If $q \geq 3$, then $p \geq 5$ is impossible. It then follows easily that $p = 3, q = 2$. Similarly, solving the equation $rs(r^2 - s^2) = 60$ in positive integers, we see that $r, s$ are divisors of 60. Say $r$ is larger, and $r + s > r > s > r - s > 0$ or $r + s > r > r - s > s > 0$; if $r \geq 5$, then $r + s \geq 6$, then this is impossible . Thus $r \leq 4$, and we have the solution $r = 4, s = 1$. Up to order and sign, then, these parameters describe the unique positive solution of the theorem.

For the second smallest solution we solve $pq(p^2-q^2)\lambda^2 = 60$ as above to find $p = 4, q = 1, \lambda = 1$. However, to solve $rs(r^2 - s^2)\mu^2 = 120$. If $\mu = 1$ we arrange so that $r + s > r > s > r - s > 0$ or $r + s > r > r - s > s > 0$, and easily find that $r = 5, s = 1$ is the only solution. If $\mu = 2$ then as in the previous case, we find $r = 3, s = 2$. These parameters give the solution stated. ∎

Next we describe a method to produce an infinite number of different integer solutions of perfect hexahedra. We consider the 'primitve'

equation, where $\lambda = \mu = 1$,
$$2(p^2 - q^2)pq = (r^2 - s^2)rs,$$
and look for a curve lying on the surface, expressed in terms of the parameters of $\alpha, \beta$; for example, in the $(p, r)$ directions, this would mean, $2(\alpha^2 - q^2)\alpha q = (\beta^2 - s^2)\beta s$.

Here is one way to do that. Suppose that $(p_0, q_0, r_0, s_0)$ is a rational solution; then $(q, q_0, s, s_0)$ for any $(q, s) \in \{(\pm q_0, 0), (0, 0), (0, \pm s_0), (\pm q_0, \pm s_0)\}$ is also a solution. Given any two rational solutions in this set of nine, express the line passing through them as an equation in terms of $x, y$. The solutions for $(x, y)$ meeting the surface, gives new solutions, $(x, q_0, y, s_0)$. Of these 9 known points, lines may pass through three of these points, and then meet the surface again at its points at infinity. However, several lines meet these 9 at only two points; these give rise to new rational solutions. For example, there are lines from the 'first quadrant': the line through $(-q_0, 0), (q_0, s_0)$, or $(0, -s_0), (q_0, s_0)$, or $(0, -s_0), (q_0, 0)$, or $(-q_0, -s_0), (q_0, 0)$. In the first case, the line is $y = \frac{s_0}{2}(\frac{x}{q_0}+1)$, which meets the surface when $x = \frac{3q_0 s_0^4}{16q_0^4 - s_0^4}$, $y = \frac{s_0(s_0^4 + 8q_0^4)}{16q_0^4 - s_0^4}$.

Now, we change $q_0, s_0$ to the parameters $q, s$ and clear denominators to obtain a parameterized curve on the surface. In a similar way, we obtain the five other parameterized curves.

(1) $\qquad (3qs^4, q(16q^4 - s^4), s(s^4 + 8q^4), s(16q^4 - s^4))$

(2) $\qquad (-3qs^4, q(16q^4 + s^4), -s(-s^4 + 8q^4), s(16q^4 + s^4))$

(3) $\qquad (-q(q^4 + 2s^4), q(q^4 - 4s^4), -3qs^4, s(q^4 - 4s^4))$

(4) $\qquad (q(q^4 - 2s^4), q(q^4 + 4s^4), -3qs^4, s(q^4 + 4s^4))$

(5) $\qquad (-2q(s^4 - q^4), q(2q^4 + s^4), s(-s^4 + 4q^4), s(2q^4 + s^4))$

(6) $\qquad (2q(s^4 + q^4), q(2q^4 - s^4), s(s^4 + 4q^4), s(2q^4 - s^4)).$

These six different rational quintic curves, can be slightly modified with sign changes; however the sign changes yield essentially the same solutions for $a, b, c, d, e, f$. These may be the smallest degree rational curves on the surface which also contain rational points. This could account for the sparsity of solutions as observed in [1].

For the second parameterization above, we can view it in 3-space by dehomogenizing with respect to the last variable and letting $t = \frac{q}{s}$ to give $(\frac{-3t}{16t^4+1}, t, \frac{1-8t^4}{16t^4+1})$. In fact by projection into the plane of the first $x$ and third $z$ variable, we obtain the curve with equation $8x^4 + 8z^4 + 4z^3 - 6z^2 - 5z - 1 = 0$, [Figure 2], which however when in 3-space is unbounded but having asymptote $x = 0, z = \frac{-1}{2}$, [Figure 1]. The quartic surface is also displayed. It has an easily discernible hole; the curve follows a fold and over the central bridge, then along an opposite fold in the surface.

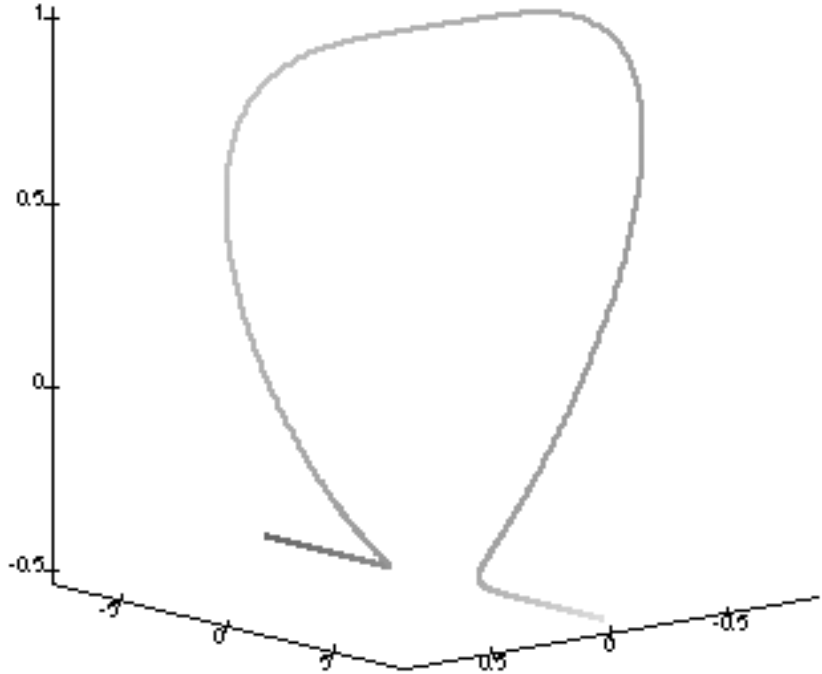

FIGURE 1. Parameterized Curve of Solutions in 3-space

We obtain integer solutions for $a, b, c, d, e, f$ using the parametrizations discussed above; for an integer $q = n$, $s = 1$ and using this particular curve discussed we have the positive solutions

$$a = 48n^4(4n^4 + 1), b = 2(-1 + 8n^4)(1 + 16n^4), c = 2 + 16n^4 + 320n^8$$

and

$$e = 2n^2(-7 - 32n^4 + 128n^8), f = 2n^2(-1 + 64n^4 + 128n^8), d = 2n^2(5 + 16n^4 + 128n^8)$$

which give infinitely many distinct integral solutions on the curve since $\frac{a}{b}, \frac{e}{f}$ yield infinitely many distinct rationals.

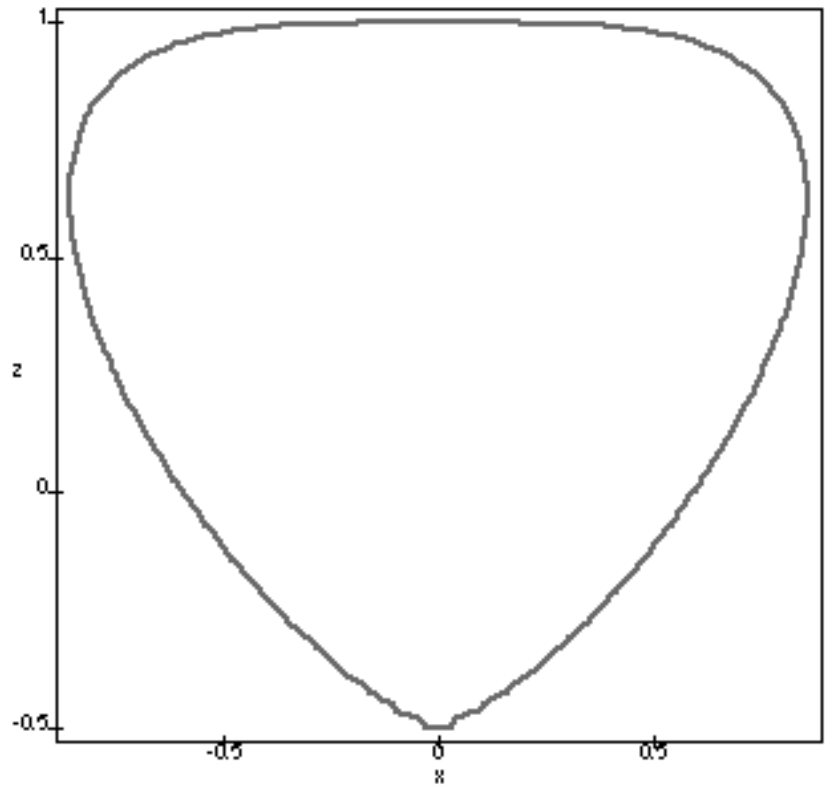

Figure 2. Projection of Space Curve into x-z Plane

An easy calculation of the surface area for these examples show that it grows like a polynomial function in $n^{18}$. Also, we can position the hexahedron in space with the center of both the top and bottom rectangles meeting the z-axis; a corner of the bottom rectangle is placed at $(-\frac{a}{2}, \frac{b}{2}, 0)$ and the top rectangle is placed with corner at $(\frac{b}{2}, \frac{a}{2}, h)$. We find that the height $h$ satisfies $a^2 + b^2 + 2h^2 = 2d^2 = e^2 + f^2$. From this we see that, $h^2 = 2(2n^4 - 1)(-1 + 8n^4)(32n^4 - 1)(64n^8 + 8n^4 + 1)$; thus the height is non-zero. Summarizing this we have the following.

**Theorem 2.** *There are infinitely many solutions to $a^2 + b^2 = c^2, d^2 = e^2 + ab, f^2 = d^2 + ab$. Thus, there are infinitely many non-trivial dissimilar perfect hexahedra with two rectangular congruent opposite parallel faces, and four congruent trapezoidal faces.*

E-mail: alperin@mathcs.sjsu.edu, Department of Mathematics and Computer Science, San Jose State University, San Jose, CA 95192 USA

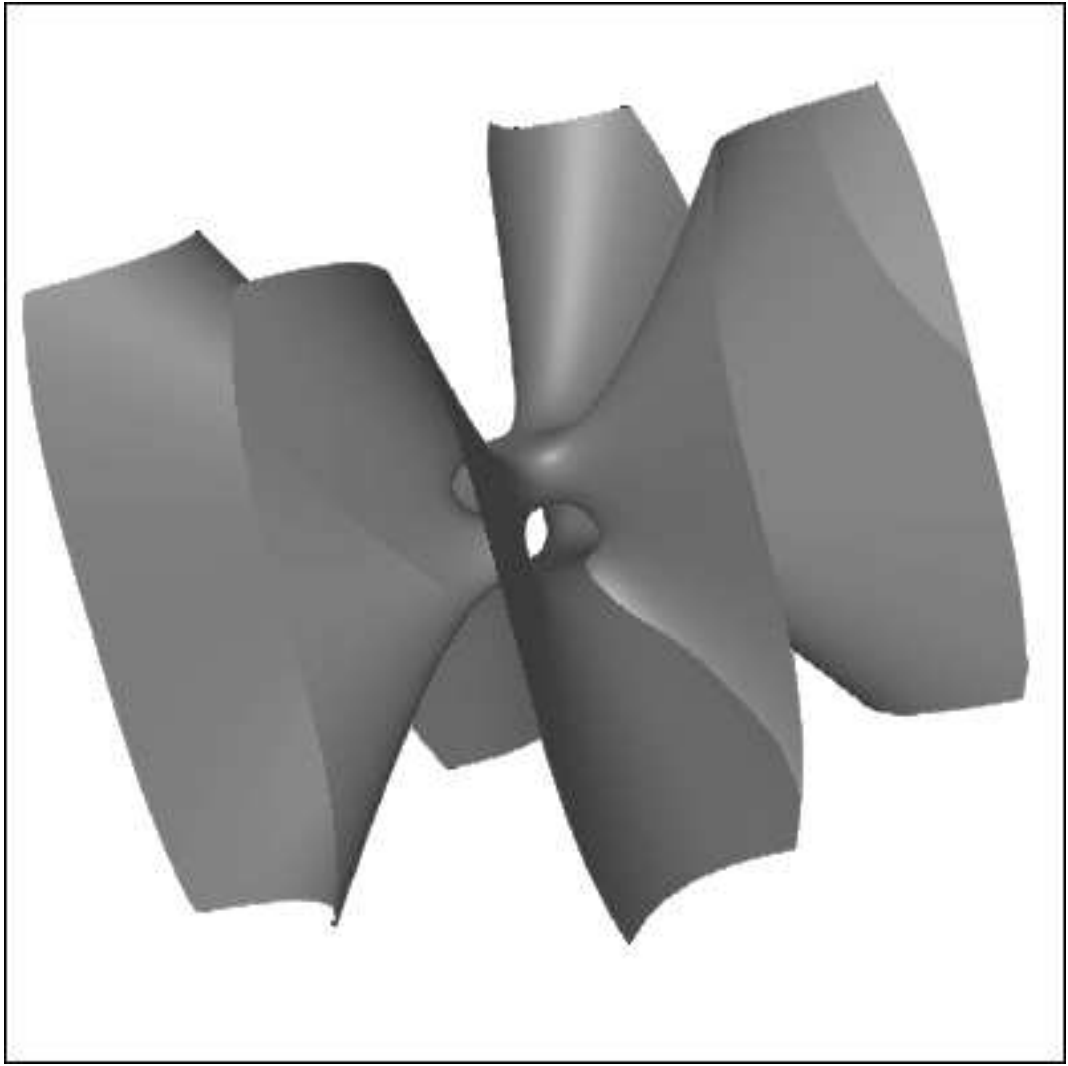

Figure 3. Quartic Surface